\documentclass[12pt]{article}
\usepackage{amssymb, amsmath}

\newtheorem{definition}{Definition}[section]
\newtheorem{theorem}[definition]{Theorem}

\newtheorem{lemma}[definition]{Lemma}
\newtheorem{proposition}[definition]{Proposition}

\newtheorem{remark}[definition]{Remark}
\newtheorem{corollary}[definition]{Corollary}

\newtheorem{notation}[definition]{Notation}
\newtheorem{assumption}[definition]{Assumption}
\def\R{\mathbb R}
\def\C{\mathbb C}
\def\Z{\mathbb Z}

\newcommand{\beast}{\begin{eqnarray*}}
\newcommand{\eeast}{\end{eqnarray*}}

\begin{document}

\title{\bf Some matrices associated with the split decomposition for a $Q$-polynomial distance-regular graph}

\author{
Joohyung Kim
 }

\date{}

\maketitle

\begin{abstract}
We consider a $Q$-polynomial distance-regular graph $\Gamma$ with
vertex set $X$ and diameter $D \geq 3$. For $\mu, \nu \in \lbrace
\downarrow, \uparrow \rbrace$ we define a direct sum decomposition
of the standard module $V=\C X$, called the $(\mu,\nu)$--split
decomposition. For this decomposition we compute the complex
conjugate and transpose of the associated primitive idempotents. Now
fix $b,\beta \in \mathbb C$ such that $b \neq 1$ and assume $\Gamma$
has classical parameters $(D,b,\alpha,\beta)$ with $\alpha = b-1$.
Under this assumption Ito and Terwilliger displayed an action of the
$q$-tetrahedron algebra $\boxtimes_q$ on the standard module of
$\Gamma$. To describe this action they defined eight matrices in
 $\hbox{Mat}_X(\mathbb C)$, called
\begin{eqnarray*}
\label{eq:list} A,\quad A^*,\quad B,\quad B^*, \quad K,\quad
K^*,\quad \Phi,\quad \Psi.
\end{eqnarray*}
For each matrix in the above list we compute the transpose and
complex conjugate. Using this information we compute the transpose
and complex conjugate for each generator of $\boxtimes_q$ on $V$.

\medskip
\noindent {\bf Keywords}. Distance-regular graph, Terwilliger
algebra, split decomposition, $q$-tetrahedron algebra.
 \hfil\break
\noindent {\bf 2000 Mathematics Subject Classification}. Primary
05E30; Secondary 05E35, 05C50
\end{abstract}

\section{Introduction}

We consider a $Q$-polynomial distance-regular graph $\Gamma$ with
vertex set $X$ and diameter $D \geq 3$ (see Section 2 for formal
definitions). Let $V$ denote the vector space over $\C$ consisting
of column vectors whose coordinates are indexed by $X$ and whose
entries are in $\C$. We call $V$ the {\it standard module}. For
$\mu, \nu \in \lbrace \downarrow, \uparrow \rbrace$ we define a
direct sum decomposition of $V$ called the $(\mu,\nu)$--{\it split
decomposition}. In this decomposition the components are indexed by
ordered pairs of integers ($i,j$) ($0 \leq i,j\leq D$). For $0 \leq
i,j\leq D$ we let $E^{\mu \nu}_{i,j}$ denote the projection onto
component ($i,j$) of the decomposition. We show that $E^{\mu
\nu}_{i,j}$ is real. We show that the transpose $(E^{\downarrow
\downarrow}_{i,j})^t=E^{\uparrow \uparrow}_{D-i,D-j}$ and
$(E^{\downarrow \uparrow}_{i,j})^t=E^{\uparrow
\downarrow}_{D-i,D-j}$. Now fix $b,\beta \in \mathbb C$ such that
$b\neq 1$ and assume $\Gamma$ has classical parameters
$(D,b,\alpha,\beta)$ with $\alpha=b-1$. Fix $q \in \C$ such that $b
= q^2$. Under this assumption Ito and Terwilliger displayed an
action of the $q$-tetrahedron algebra $\boxtimes_q$ on the standard
module of $\Gamma$ \cite{drgqtet}. To describe this action they
defined eight matrices in
 $\hbox{Mat}_X(\mathbb C)$, called
\begin{eqnarray*}
\label{eq:list} A,\quad A^*,\quad B,\quad B^*, \quad K,\quad
K^*,\quad \Phi,\quad \Psi.
\end{eqnarray*}
For each matrix in the above list we compute the transpose and
complex conjugate. Concerning the transpose, we show that each of
$A$, $A^*$, $\Phi$, $\Psi$ is symmetric and that $B^t=B^*$ and
$K^t={K^*}^{-1}$. Concerning the complex conjugate, recall from the
above assumption that $q \in \C$ satisfies $b = q^2$. Define
$q{'}=-q$ and note that $(q')^2=b$. There are two cases to consider.
For $b>1$ we have $\overline{q}=q$ and $\overline{q'}=q'$. For
$b<-1$ we have $q \in i\mathbb R$ so $\overline{q}=q'$. For each
matrix $S$ in the above list let $S'$ denote the corresponding
matrix associated with $q'$. We show that for each matrix $S$ from
the above list, $S$ is real if $b>1$ and $\overline{S} = S'$ if
$b<-1$. Using the above information we compute the transpose and
complex conjugate for each generator of $\boxtimes_q$ on $V$.

\section{The Terwilliger algebra of a distance-regular graph; preliminaries}
\noindent In this section we review some definitions and basic
concepts concerning the Terwilliger algebra of a distance-regular
graph. For more background information we refer the reader to
\cite{bannai, bcn, godsil, terwSub1}.

\medskip
\noindent Let $X$ denote a nonempty  finite  set. Let
 $\hbox{Mat}_X(\C)$
denote the $\C$-algebra consisting of all matrices whose rows and
columns are indexed by $X$ and whose entries are in $\C$. Let $V=\C
X$ denote the vector space over $\C$ consisting of column vectors
whose coordinates are indexed by $X$ and whose entries are in $\C$.
We observe that $\hbox{Mat}_X(\C)$ acts on $V$ by left
multiplication. We call $V$ the {\it standard module}. We endow $V$
with the Hermitian form $\langle \, , \, \rangle$ that satisfies
$\langle u,v \rangle = u^t\overline{v}$ for $u,v \in V$, where $t$
denotes transpose. Observe that $\langle \, , \, \rangle$ is
positive definite. We call this form the {\it standard Hermitian
form} on $V$. Observe that for $B \in \hbox{Mat}_X(\C)$,
\begin{equation}
\label{SHF} \langle Bu, v\rangle =\langle u, {\overline{B}}^{t}v
\rangle \qquad \qquad  u,v \in V.
\end{equation}
For all $y \in X,$ let $\hat{y}$ denote the element of $V$ with a 1
in the $y$ coordinate and 0 in all other coordinates. Observe that
$\{\hat{y}\;|\;y \in X\}$ is an orthonormal basis for $V.$

\medskip
\noindent Let $\Gamma = (X,R)$ denote a finite, undirected,
connected graph, without loops or multiple edges, with vertex set
$X$ and edge set $R$. Let $\partial $ denote the path-length
distance function for $\Gamma $,  and set $D :=
\mbox{max}\{\partial(x,y) \;|\; x,y \in X\}$. We call $D$  the {\it
diameter} of $\Gamma $. We say that $\Gamma$ is {\it
distance-regular} whenever for all integers $h,i,j\;(0 \le h,i,j \le
D)$ and for all vertices $x,y \in X$ with $\partial(x,y)=h,$ the
number
\begin{eqnarray*}
p_{ij}^h = |\{z \in X \; |\; \partial(x,z)=i, \partial(z,y)=j \}|
\end{eqnarray*}
is independent of $x$ and $y.$ The $p_{ij}^h$ are called the {\it
intersection numbers} of $\Gamma.$ We abbreviate $c_i=p^i_{1,i-1}$
$(1 \leq i \leq D)$, $b_i=p^i_{1,i+1}$ $(0 \leq i \leq D-1)$,
$a_i=p^i_{1i}$ $(0 \leq i \leq D)$.

\medskip
\noindent For the rest of this paper we assume  that $\Gamma$ is
distance-regular  with diameter $D\geq 3$.

\medskip
\noindent We recall the Bose-Mesner algebra of $\Gamma.$ For $0 \le
i \le D$ let $A_i$ denote the matrix in $\hbox{Mat}_X(\C)$ with $xy$
entry
$$
(A_i)_{xy} = \begin{cases}
                  1, & \text{if } \partial(x,y)=i \\
                  0, & \text{if } \partial(x,y) \ne i
             \end{cases}
 \qquad (x,y \in X).
$$
We call $A_i$ the $i$th {\it distance matrix} of $\Gamma.$ $A_1$ is
called the {\it adjacency matrix} of $\Gamma.$ We observe that (i)
$A_0 = I$;
 (ii) $\sum_{i=0}^D A_i = J$;
(iii) $\overline{A_i} = A_i \;(0 \le i \le D)$; (iv) $A_i^t = A_i
\;(0 \le i \le D)$; (v) $A_iA_j = \sum_{h=0}^D p_{ij}^h A_h \;( 0
\le i,j \le D)$, where $I$ (resp. $J$) denotes the identity matrix
(resp. all 1's matrix) in
 $\hbox{Mat}_X(\C)$.
 Using these facts  we find $A_0,A_1,\ldots,A_D$
form a basis for a commutative subalgebra $M$ of $\mbox{Mat}_X(\C)$.
We call $M$ the {\it Bose-Mesner algebra} of $\Gamma$. It turns out
that $A_1$ generates $M$ \cite[p.~190]{bannai}. By
\cite[p.~45]{bcn}, $M$ has a second basis $E_0,E_1,\ldots,E_D$ such
that (i) $E_0 = |X|^{-1}J$; (ii) $\sum_{i=0}^D E_i = I$; (iii)
$\overline{E_i} = E_i \;(0 \le i \le D)$; (iv) $E_i^t =E_i  \;(0 \le
i \le D)$; (v) $E_iE_j =\delta_{ij}E_i  \;(0 \le i,j \le D)$. We
call $E_0, E_1, \ldots, E_D $  the {\it primitive idempotents} of
$\Gamma$.

\medskip
\noindent We recall the eigenvalues of $\Gamma $. Since
$E_0,E_1,\ldots,E_D$ form a basis for $M$ there exist complex
scalars $\theta_0,\theta_1, \ldots,\theta_D$ such that $A_1 =
\sum_{i=0}^D \theta_iE_i$. Observe that $A_1 E_i = E_i A_1 =
\theta_iE_i$ for $0 \leq i \leq D$. We call $\theta_i$  the {\it
eigenvalue} of $\Gamma$ associated with $E_i$ $(0 \leq i \leq D)$.
Since $A_1$ is real symmetric, the eigenvalues
$\theta_0,\theta_1,\ldots,\theta_D$ are in $\R.$ Observe that
$\theta_0,\theta_1,\ldots,\theta_D$ are mutually distinct since
$A_1$ generates $M$. Observe that
\begin{eqnarray}
\label{vsub} V = E_0V+E_1V+ \cdots +E_DV \qquad \qquad {\rm
(orthogonal\ direct\ sum}).
\end{eqnarray}
For $0 \le i \le D$ the space $E_iV$ is the  eigenspace of $A_1$
associated with $\theta_i$.

\medskip
\noindent We recall the Krein parameters. Let $\circ $ denote the
entrywise product in $\mbox{Mat}_X(\C)$. Observe that $A_i\circ A_j=
\delta_{ij}A_i$ for $0 \leq i,j\leq D$, so $M$ is closed under
$\circ$. Thus there exist complex scalars $q^h_{ij}$ $(0 \leq
h,i,j\leq D)$ such that
$$
E_i\circ E_j = |X|^{-1}\sum_{h=0}^D q^h_{ij}E_h \qquad (0 \leq
i,j\leq D).
$$
By \cite[p.~170]{Biggs}, $q^h_{ij}$ is real and nonnegative  for $0
\leq h,i,j\leq D$. The $q^h_{ij}$ are called the {\it Krein
parameters}. The graph $\Gamma$ is said to be {\it $Q$-polynomial}
(with respect to the given ordering $E_0, E_1, \ldots, E_D$ of the
primitive idempotents) whenever for $0 \leq h,i,j\leq D$, $q^h_{ij}=
0$ (resp. $q^h_{ij}\not= 0$) whenever one of $h,i,j$ is greater than
(resp. equal to) the sum of the other two \cite[p.~59]{bcn}. See
\cite{bannai, caugh1, caugh2, aap1, aap2} for more information on
the $Q$-polynomial property. From now on we assume that $\Gamma$ is
$Q$-polynomial with respect to $E_0, E_1, \ldots, E_D$.

\medskip
\noindent We recall the dual Bose-Mesner algebra of $\Gamma.$ Fix a
vertex $x \in X.$ We view $x$ as a ``base vertex''. For $ 0 \le i
\le D$ let $E_i^*=E_i^*(x)$ denote the diagonal matrix in
$\hbox{Mat}_X(\C)$ with $yy$ entry

\begin{equation}\label{DEFDEI}
(E_i^*)_{yy} = \begin{cases}
                      1, & \text{if } \partial(x,y)=i \\
                      0, & \text{if } \partial(x,y) \ne i
                \end{cases}
 \qquad (y \in X).
\end{equation}
We call $E_i^*$ the  $i$th {\it dual idempotent} of $\Gamma$
 with respect to $x$ \cite[p.~378]{terwSub1}.
We observe that (i) $\sum_{i=0}^D E_i^*=I$; (ii) $\overline{E_i^*} =
E_i^*$ $(0 \le i \le D)$; (iii) $E_i^{*t} = E_i^*$ $(0 \le i \le
D)$; (iv) $E_i^*E_j^* = \delta_{ij}E_i^* $ $(0 \le i,j \le D)$. By
these facts $E_0^*,E_1^*, \ldots, E_D^*$ form a basis for a
commutative subalgebra $M^*=M^*(x)$ of $\hbox{Mat}_X(\C).$ We call
$M^*$ the {\it dual Bose-Mesner algebra} of $\Gamma$ with respect to
$x$ \cite[p.~378]{terwSub1}. For $0 \leq i \leq D$ let $A^*_i =
A^*_i(x)$ denote the diagonal matrix in
 $\hbox{Mat}_X(\C)$
with $yy$ entry $(A^*_i)_{yy}=\vert X \vert (E_i)_{xy}$ for $y \in
X$. Then $A^*_0, A^*_1, \ldots, A^*_D$ form a basis for $M^*$
\cite[p.~379]{terwSub1}. Moreover (i) $A^*_0 = I$; (ii)
$\overline{A^*_i} = A^*_i \;(0 \le i \le D)$; (iii) $A^{*t}_i =
A^*_i  \;(0 \le i \le D)$; (iv) $A^*_iA^*_j = \sum_{h=0}^D q_{ij}^h
A^*_h \;( 0 \le i,j \le D) $ \cite[p.~379]{terwSub1}. We call
 $A^*_0, A^*_1, \ldots, A^*_D$
the {\it dual distance matrices} of $\Gamma$ with respect to $x$.
$A^*_1$ is called the {\it dual adjacency matrix} of $\Gamma$ with
respect to $x$. The matrix $A^*_1$ generates $M^*$
\cite[Lemma~3.11]{terwSub1}.

\medskip
\noindent We recall the dual eigenvalues of $\Gamma$. Since
$E^*_0,E^*_1,\ldots,E^*_D$ form a basis for $M^*$ and since $A^*_1$
is real, there exist real scalars $\theta^*_0,\theta^*_1,
\ldots,\theta^*_D$ such that $A^*_1 = \sum_{i=0}^D \theta^*_iE^*_i$.
Observe that $A^*_1 E^*_i = E^*_iA^*_1 =  \theta^*_iE^*_i$ for $0
\leq i \leq D$.
 We call $\theta^*_i$ the {\it dual eigenvalue}
of $\Gamma$ associated with $E^*_i$ $(0 \leq i\leq D)$. Observe that
$\theta^*_0,\theta^*_1,\ldots,\theta^*_D$ are mutually distinct
since $A^*_1$ generates $M^*$.

\medskip
\noindent We recall the subconstituents of $\Gamma $. From
(\ref{DEFDEI}) we find
\begin{equation}\label{DEIV}
E_i^*V = \mbox{span}\{\hat{y} \;|\; y \in X, \quad \partial(x,y)=i\}
\qquad (0 \le i \le D).
\end{equation}
By (\ref{DEIV})  and since
 $\{\hat{y}\;|\;y \in X\}$ is an orthonormal basis for $V$
 we find
\begin{eqnarray}
\label{vsubstar} V = E_0^*V+E_1^*V+ \cdots +E_D^*V \qquad \qquad
{\rm (orthogonal\ direct\ sum}).
\end{eqnarray}
For $0 \leq i \leq D$ the space $E^*_iV$ is the eigenspace of
$A^*_1$ associated with $\theta^*_i$. We call $E_i^*V$ the $i$th
{\it subconstituent} of $\Gamma$ with respect to $x$.

\medskip
\noindent We recall the subconstituent algebra of $\Gamma $. Let
$T=T(x)$ denote the subalgebra of $\hbox{Mat}_X(\C)$ generated by
$M$ and $M^*$. We call $T$ the {\it subconstituent algebra} (or {\it
Terwilliger algebra}) of $\Gamma$
 with respect to $x$ \cite[Definition~3.3]{terwSub1}.
We observe that $T$ is generated by $A_1,A^*_1$. We observe that $T$
has finite dimension. Moreover $T$ is semi-simple since it is closed
under the conjugate transpose map \cite[p.~157]{CR}. See
\cite{curtin1, go2, terwSub1, terwSub2, terwSub3} for more
information on the subconstituent algebra.

\medskip
\noindent Until further notice we adopt the following notation.

\begin{notation}
\label{setup} \rm We assume that $\Gamma=(X,R)$ is a
distance-regular graph with diameter $D\geq 3$. We assume that
$\Gamma$ is $Q$-polynomial with respect to the ordering $E_0, E_1,
\ldots, E_D$ of the primitive idempotents. We fix $x \in X$ and
write $A^*_1=A^*_1 (x)$, $E^*_i=E^*_i(x)$ $(0 \leq i \leq D)$,
$T=T(x)$. We abbreviate $V=\C X$.
\end{notation}

\noindent We recall some useful results on $T$-modules. With
reference to Notation \ref{setup}, by a {\it T-module} we mean a
subspace $W \subseteq V$ such that $BW \subseteq W$ for all $B \in
T.$
 Let $W$ denote a $T$-module. Then $W$ is said
to be {\it irreducible} whenever $W$ is nonzero and $W$ contains no
$T$-modules other than 0 and $W$. Let $W$ denote a $T$-module and
let $W'$ denote a $T$-module contained in $W$. Then the orthogonal
complement of $W'$ in $W$ is a $T$-module \cite[p.~802]{go2}. It
follows that each $T$-module is an orthogonal direct sum of
irreducible $T$-modules. In particular $V$ is an orthogonal direct
sum of irreducible $T$-modules. Let $W$ denote an irreducible
$T$-module. By the {\it endpoint} of $W$ we mean $\mbox{min}\lbrace
i |0\leq i \leq D, \; E^*_iW\not=0\rbrace $. By the {\it diameter}
of $W$ we mean $ |\lbrace i | 0 \leq i \leq D,\; E^*_iW\not=0
\rbrace |-1 $. By the {\it dual endpoint} of $W$ we mean
$\mbox{min}\lbrace i |0\leq i \leq D, \; E_iW\not=0\rbrace $. By the
{\it dual diameter} of $W$ we mean $ |\lbrace i | 0 \leq i \leq D,\;
E_iW\not=0 \rbrace |-1 $. The diameter of $W$ is equal to the dual
diameter of $W$ \cite[Corollary~3.3]{aap1}.

\begin{lemma}
{\rm \cite[Lemma~3.4, Lemma~3.9, Lemma~3.12]{terwSub1}}
\label{lem:basic} With reference to Notation \ref{setup}, let $W$
denote an irreducible $T$-module with endpoint $\rho$, dual endpoint
$\tau$, and diameter $d$. Then $\rho,\tau,d$ are nonnegative
integers such that $\rho+d\leq D$ and $\tau+d\leq D$. Moreover the
following (i)--(iv) hold.
\begin{enumerate}
\item
$E^*_iW \not=0$ if and only if $\rho \leq i \leq \rho+d$ $ \quad (0
\leq i \leq D)$.
\item
$W = \sum_{h=0}^{d} E^*_{\rho+h}W$ \qquad \rm({orthogonal direct
sum}).
\item
$E_iW \not=0$ if and only if $\tau \leq i \leq \tau+d$ $ \quad (0
\leq i \leq D)$.
\item
$W = \sum_{h=0}^{d} E_{\tau+h}W$ \qquad \rm({orthogonal direct
sum}).
\end{enumerate}
\end{lemma}

\begin{lemma}
{\rm \cite[Lemma 5.1, Lemma 7.1]{caugh2}} \label{lem:basineq} With
reference to Notation \ref{setup}, let $W$ denote an irreducible
$T$-module with endpoint $\rho$, dual endpoint $\tau$, and diameter
$d$. Then the following (i), (ii) hold.
\begin{enumerate}
\item $2\rho+d\geq D$. \item $2\tau+d\geq D$.
\end{enumerate}
\end{lemma}

\medskip
\noindent We finish this section with a comment.

\begin{lemma}
{\rm \cite[Lemma~3.3]{duality}} \label{lem:wsplit} With reference to
Notation \ref{setup}, let $W$ denote an irreducible $T$-module with
endpoint $\rho$, dual endpoint $\tau$, and diameter $d$. Then for
$\mu, \nu \in \lbrace \downarrow,\uparrow \rbrace$ we have
\begin{eqnarray*}
\label{eq:splitd} W = \sum_{h=0}^d W^{\mu \nu}_h \qquad (\rm{direct\
sum}),
\end{eqnarray*}
where for $0 \leq h\leq d$,
\begin{eqnarray*}
W^{\downarrow \downarrow}_h &=& (E^*_\rho W+\cdots+
E^*_{\rho+h}W)\cap (E_\tau W+\cdots + E_{\tau+d-h}W),
\label{eq:whdef1}
\\
W^{\uparrow \downarrow}_h &=& (E^*_{\rho+d-h} W+\cdots+
E^*_{\rho+d}W)\cap (E_\tau W+\cdots + E_{\tau+d-h}W),
\label{eq:whdef2}
\\
W^{\downarrow \uparrow}_h &=& (E^*_\rho W+\cdots+ E^*_{\rho+h}W)\cap
(E_{\tau+h}W+\cdots + E_{\tau+d}W), \label{eq:whdef3}
\\
W^{\uparrow \uparrow}_h &=& (E^*_{\rho+d-h} W+\cdots+
E^*_{\rho+d}W)\cap (E_{\tau+h}W+\cdots + E_{\tau+d}W).
\label{eq:whdef4}
\end{eqnarray*}
\end{lemma}

\section{The displacement and split decompositions of the standard module}

\noindent In this section we recall the displacement decompositions
and the four split decompositions for the standard module and
discuss their basic properties and connections. We start with a
definition.

\begin{definition}
{\rm \cite[Definition~4.1]{ds}} \label{def:shift} \rm With reference
to Notation
 \ref{setup}, let $W$ denote an irreducible $T$-module.
By the {\it displacement of $W$ of the first kind} (resp. {\it
second kind}) we mean the integer $\rho+\tau+d-D$ (resp.
$\rho-\tau$) where $\rho,\tau,d$ denote respectively the endpoint,
dual endpoint, and diameter of $W$.
\end{definition}

\begin{lemma}
\label{lem:shiftpos} With reference to Notation \ref{setup}, let $W$
denote an irreducible $T$-module. Then the following (i), (ii) hold.
\begin{enumerate}
\item Let $\eta$ denote the displacement of $W$ of the first kind. Then $0 \leq \eta \leq D$.
\item Let $\zeta$ denote the displacement of $W$ of the second kind. Then $-D \leq \zeta \leq D$.
\end{enumerate}
\end{lemma}

\noindent {\it Proof:} (i) This is just \cite[Lemma~4.2]{ds}.
\\
(ii) Recall that $0 \leq \rho,\tau \leq D$ by Lemma \ref{lem:basic}.
The result follows.
\hfill $\Box$\\

\begin{definition}
\label{def:modl1} \rm With reference to Notation \ref{setup}, for $0
\leq \eta \leq D$ let ${V}_\eta$ denote the subspace of $V$ spanned
by the irreducible $T$-modules for which $\eta$ is the displacement
of the first kind. Observe that $V_\eta$ is a $T$-module. By
\cite[Lemma~4.4]{ds} we have $V = \sum_{\eta=0}^D {V}_\eta
~(\rm{orthogonal\ direct\ sum})$. For $0 \leq \eta \leq D$ we define
a matrix $\phi_{\eta} \in \hbox{Mat}_X(\C)$ so that
\begin{eqnarray*}
\label{phiV1} (\phi_{\eta} - I)V_{\eta}&=&0,
\\
\label{phiV2} \phi_{\eta}V_{\xi}&=&0 \quad \hbox{if}~\eta \neq \xi
\qquad (0 \leq \xi \leq D).
\end{eqnarray*}
In other words $\phi_{\eta}$ is the projection from $V$ onto
${V}_\eta$. We note that $V_{\eta} = \phi_{\eta} V$.
\end{definition}

\medskip
\noindent The following result is immediate from Definition
\ref{def:modl1}.

\begin{lemma}
\label{etads}
 With reference to Notation \ref{setup},
\begin{eqnarray}
\label{V-eta} V = \sum_{\eta=0}^D {\phi}_\eta V \qquad \qquad
(\rm{orthogonal\ direct\ sum}).
\end{eqnarray}
Moreover for $0 \leq \eta \leq D$ ${\phi}_\eta V$ is the subspace of
$V$ spanned by the irreducible $T$-modules for which $\eta$ is the
displacement of the first kind.
\end{lemma}

\begin{definition}
\rm We call the sum (\ref{V-eta}) the {\it displacement
decomposition of $V$ of the first kind} with respect to $x$.
\end{definition}

\begin{lemma}
\label{lem:prosumphiV1} With reference to Notation \ref{setup} and
Definition \ref{def:modl1}, the following (i), (ii) hold.
\begin{enumerate}
\item $\sum_{\eta=0}^D  \phi_{\eta}=I$.
\item $\phi_{\eta}\phi_{\xi}=\delta_{\eta \xi} \phi_{\eta}$  \quad {\rm ($0 \leq \eta,\xi \leq
D$)}.
\end{enumerate}
\end{lemma}

\noindent {\it Proof:} Immediate from Definition \ref{def:modl1}.
\hfill $\Box$\\

\medskip
\noindent We have been discussing the displacement of the first
kind. We now do something similar for the displacement of the second
kind.

\begin{definition}
\label{def:modl2} \rm With reference to Notation \ref{setup}, for
$-D \leq \zeta \leq D$ let ${V}_\zeta$ denote the subspace of $V$
spanned by the irreducible $T$-modules for which $\zeta$ is the
displacement of the second kind. Observe that $V_\zeta$ is a
$T$-module. By \cite[Lemma~4.4]{ds} $V = \sum_{\zeta=-D}^D {V}_\zeta
~(\rm{orthogonal\ direct\ sum})$. For $-D \leq \zeta \leq D$ we
define a matrix $\psi_{\zeta} \in \hbox{Mat}_X(\C)$ so that
\begin{eqnarray*}
\label{psiV1} (\psi_{\zeta} - I)V_{\zeta}&=&0,
\\
\label{psiV2} \psi_{\zeta}V_{\xi}&=&0 \quad \hbox{if}~\zeta \neq \xi
\qquad (-D \leq \xi \leq D).
\end{eqnarray*}
In other words $\psi_{\zeta}$ is the projection from $V$ onto
${V}_\zeta$. We note that ${V}_\zeta = \psi_{\zeta}V.$
\end{definition}

\medskip
\noindent The following result is immediate from Definition
\ref{def:modl2}.

\begin{lemma}
\label{zetads}
 With reference to Notation \ref{setup},
\begin{eqnarray}
\label{V-zeta} V = \sum_{\zeta=-D}^D {\psi}_\zeta V \qquad \qquad
(\rm{orthogonal\ direct\ sum}).
\end{eqnarray}
Moreover for $-D \leq \zeta \leq D$ ${\psi}_\zeta V$ is the subspace
of $V$ spanned by the irreducible $T$-modules for which $\zeta$ is
the displacement of the second kind.
\end{lemma}

\begin{definition}
\rm We call the sum (\ref{V-zeta}) the {\it displacement
decomposition of $V$ of the second kind} with respect to $x$.
\end{definition}

\begin{lemma}
\label{lem:prosumpsiV1} With reference to Notation \ref{setup} and
Definition \ref{def:modl2}, the following (i), (ii) hold.
\begin{enumerate}
\item $\sum_{\zeta={-D}}^D  \psi_{\zeta}=I$.
\item $\psi_{\zeta}\psi_{\xi}=\delta_{\zeta \xi} \psi_{\zeta}$ \quad {\rm ($-D \leq \zeta,\xi \leq
D$)}.
\end{enumerate}
\end{lemma}

\noindent {\it Proof:} Immediate from Definition \ref{def:modl2}.
\hfill $\Box$\\

\medskip
\noindent We now recall the split decompositions of $V$.

\begin{definition}
{\rm \cite[Definition~10.1]{drgqtet}} \label{def:updown} \rm With
reference to Notation \ref{setup}, for $-1\leq i,j\leq D$ we define
\begin{eqnarray*}
V^{\downarrow \downarrow}_{i,j} &=& (E^*_0V+\cdots+E^*_iV)\cap
(E_0V+\cdots+E_jV),
\\
V^{\uparrow \downarrow}_{i,j} &=& (E^*_DV+\cdots+E^*_{D-i}V)\cap
(E_0V+\cdots+E_jV),
\\
V^{\downarrow \uparrow}_{i,j} &=& (E^*_0V+\cdots+E^*_iV)\cap
(E_DV+\cdots+E_{D-j}V),
\\
V^{\uparrow \uparrow}_{i,j} &=& (E^*_DV+\cdots+E^*_{D-i}V)\cap
(E_DV+\cdots+E_{D-j}V).
\end{eqnarray*}
In each of the above four equations we interpret the right-hand side
to be 0 if $i=-1$ or $j=-1$.
\end{definition}

\begin{lemma}
\label{lem:updown1} With reference to Notation \ref{setup}, the
following hold for $0 \leq i \leq D$.
\begin{eqnarray}
\label{updown1}
 V^{\downarrow \downarrow}_{i,D} &=&
E^*_0V+\cdots+E^*_iV, \quad \quad ~~\,\,\,V^{\downarrow
\downarrow}_{D,i} ~\,=~\,
 E_0V+\cdots+E_iV,
\\
\label{updown2}
 V^{\uparrow \downarrow}_{i,D} &=&
 E^*_DV+\cdots+E^*_{D-i}V,
\quad \quad V^{\uparrow \downarrow}_{D,i} ~\,=~\,
 E_0V+\cdots+E_iV,
\\
\label{updown3}
 V^{\downarrow \uparrow}_{i,D} &=&
 E^*_0V+\cdots+E^*_iV,
\quad \quad ~~~\,V^{\downarrow \uparrow}_{D,i} ~\,=~\,
 E_DV+\cdots+E_{D-i}V,
\\
\label{updown4}
 V^{\uparrow \uparrow}_{i,D} &=&
 E^*_DV+\cdots+E^*_{D-i}V,
\quad \quad V^{\uparrow \uparrow}_{D,i} ~\,=~\,
 E_DV+\cdots+E_{D-i}V.
\end{eqnarray}
\end{lemma}

\noindent {\it Proof:} Immediate from (\ref{vsub}),
(\ref{vsubstar}), and Definition \ref{def:updown}.
 \hfill $\Box$\\

\begin{definition}
{\rm \cite[Definition~10.2]{drgqtet}} \label{def:vtilde} \rm With
reference to Notation \ref{setup} and Definition \ref{def:updown},
for $\mu, \nu \in \lbrace \downarrow, \uparrow \rbrace$ and  $0 \leq
i,j\leq D$ we have $ V^{\mu \nu}_{i-1,j} \subseteq V^{\mu
\nu}_{i,j}$ and $ V^{\mu \nu}_{i,j-1}  \subseteq V^{\mu \nu}_{i,j}
$. Therefore
\begin{eqnarray*}
V^{\mu \nu}_{i-1,j}+ V^{\mu \nu}_{i,j-1} \subseteq V^{\mu
\nu}_{i,j}.
\end{eqnarray*}
Referring to the above inclusion, we define ${\tilde V}^{\mu
\nu}_{i,j}$
 to be the orthogonal complement of the left-hand side in the
 right-hand side; that
is
\begin{eqnarray*}
 {\tilde V}^{\mu \nu}_{i,j}=(
V^{\mu \nu}_{i-1,j}+ V^{\mu \nu}_{i,j-1})^\perp \cap V^{\mu
\nu}_{i,j}.
\end{eqnarray*}
\end{definition}

\medskip
\noindent The following result is a mild generalization of
\cite[Theorem~5.7]{ds}.

\begin{proposition}
\label{prop:usplitdec} With reference to Notation \ref{setup} and
Definition \ref{def:vtilde}, for $\mu, \nu \in \lbrace \downarrow,
\uparrow \rbrace$ and $0 \leq i,j\leq D$ we have
\begin{eqnarray}
V^{\mu \nu}_{i,j} =\sum_{r=0}^i \sum_{s=0}^j {{\tilde V}^{\mu
\nu}}_{r,s} \qquad \rm({direct\ sum}). \label{eq:vsplt0}
\end{eqnarray}
\end{proposition}

\noindent {\it Proof:} For $\mu=\downarrow $, $ \nu=\downarrow $
 this is just \cite[Theorem~5.7]{ds}.
For general values of $\mu, \nu$ the argument is essentially the
same.
\hfill $\Box$\\

\begin{corollary}
\label{cor:drsum} {\rm \cite [Lemma~10.3]{drgqtet}} With reference
to Notation \ref{setup} and Definition \ref{def:vtilde}, the
following holds
 for $\mu, \nu \in \lbrace \downarrow, \uparrow \rbrace$:
\begin{eqnarray}
V = \sum_{i=0}^D \sum_{j=0}^D {\tilde V}^{\mu \nu}_{i,j} \qquad
\qquad \rm({direct\ sum}). \label{eq:vsplt}
\end{eqnarray}
\end{corollary}

\begin{definition}
\label{def:splitdv} \rm We call the sum (\ref{eq:vsplt}) the
$(\mu,\nu)$--{\it split decomposition} of $V$ with respect to $x$.
\end{definition}

\begin{remark}
\rm The decomposition (\ref{eq:vsplt}) is not orthogonal in general.
\end{remark}

\begin{corollary}
\label{cor:dim} With reference to Notation \ref{setup} and
Definition \ref{def:vtilde}, the following (i), (ii) hold.
\begin{enumerate}
\item  For $0 \leq i\leq D$ the dimension of $E^*_i V$ is equal to each of
\begin{eqnarray*}
 \sum_{j=0}^D \mbox{dim}\,{\tilde V}^{\downarrow \downarrow}_{i,j}, \qquad  \sum_{j=0}^D \mbox{dim}\,{\tilde V}^{\uparrow \downarrow}_{D-i,j},
\\
 \sum_{j=0}^D \mbox{dim}\,{\tilde V}^{\downarrow \uparrow}_{i,j}, \qquad  \sum_{j=0}^D \mbox{dim}\,{\tilde V}^{\uparrow \uparrow}_{D-i,j}.
\end{eqnarray*}

\item For $0 \leq j \leq D$ the dimension of $E_j V$ is equal to each of
\begin{eqnarray*}
 \sum_{i=0}^D \mbox{dim}\,{\tilde V}^{\downarrow \downarrow}_{i,j}, \qquad  \sum_{i=0}^D \mbox{dim}\,{\tilde V}^{\downarrow \uparrow}_{i,D-j},
\\
 \sum_{i=0}^D \mbox{dim}\,{\tilde V}^{\uparrow \downarrow}_{i,j},   \qquad  \sum_{i=0}^D \mbox{dim}\,{\tilde V}^{\uparrow \uparrow}_{i,D-j}.
\end{eqnarray*}

\end{enumerate}
\end{corollary}

\noindent {\it Proof:} To get $\mbox{dim}\, E^*_i V = \sum_{j=0}^D
\mbox{dim}\, {\tilde V}^{\downarrow \downarrow}_{i,j}$, set $j=D$ in
(\ref{eq:vsplt0}), evaluate the result by using the equation on the
left in (\ref{updown1}), and use induction on $i$. The other
equations are similarly obtained.
\hfill $\Box$\\

\begin{lemma}
{\rm \cite[Lemma~4.6]{duality}} \label{lem:DU} With reference to
Notation \ref{setup}, let $W$ denote an irreducible $T$-module with
endpoint $\rho$, dual endpoint $\tau$, and diameter $d$. Then the
following (i)--(iv) hold for $0 \leq h\leq d$ and $0 \leq i,j\leq
D$.
\begin{enumerate}
\item $W^{\downarrow \downarrow}_h\subseteq {{\tilde V}^{\downarrow \downarrow}_{i,j}}$ if and only if $i=\rho +h$ and $j=\tau+d-h$.
\item $W^{\uparrow \downarrow}_h\subseteq {{\tilde V}^{\uparrow \downarrow}_{i,j}}$ if and only if $i=D-\rho-d+h$ and $j=\tau +d-h$.
\item $W^{\downarrow \uparrow}_h\subseteq {{\tilde V}^{\downarrow \uparrow}_{i,j}}$ if and only if $i=\rho +h$ and $j=D-\tau-h$.
\item $W^{\uparrow \uparrow}_h\subseteq {{\tilde V}^{\uparrow \uparrow}_{i,j}}$ if and only if $i=D-\rho-d+h$ and $j=D-\tau-h$.
\end{enumerate}
\end{lemma}

\medskip
\noindent In the following two theorems we describe the
relationships between the displacement decompositions and the split
decompositions.

\begin{theorem}
\label{thm:connect4-1} With reference to Notation \ref{setup}, the
following (i), (ii) hold for $0 \leq \eta \leq D$.
\begin{enumerate}
\item $\phi_{\eta} V = \sum {\tilde
V}^{\downarrow \downarrow}_{i,j}$, where the sum is over all ordered
pairs $i,j$ such that $0 \leq i,j\leq D$ and $i+j=D+\eta$.
\item $\phi_{\eta} V = \sum {\tilde
V}^{\uparrow \uparrow}_{i,j}$, where the sum is over all ordered
pairs $i,j$ such that $0 \leq i,j\leq D$ and $i+j=D-\eta$.
\end{enumerate}
\end{theorem}

\noindent {\it Proof:} Part (i) is just \cite[Theorem~6.2(i)]{ds}.
Part (ii) is similarly proved.
\hfill $\Box$\\

\begin{theorem}
\label{thm:connect4-2} With reference to Notation \ref{setup}, the
following (i), (ii) hold for $-D \leq \zeta \leq D$.
\begin{enumerate}
\item
$\psi_{\zeta} V = \sum {\tilde V}^{\downarrow \uparrow}_{i,j}$,
where the sum is over all ordered pairs $i,j$ such that $0 \leq
i,j\leq D$ and $i+j=D+\zeta$.
\item $\psi_{\zeta} V = \sum {\tilde
V}^{\uparrow \downarrow}_{i,j}$, where the sum is over all ordered
pairs $i,j$ such that $0 \leq i,j\leq D$ and $i+j=D-\zeta$.
\end{enumerate}
\end{theorem}

\noindent {\it Proof:} Similar to the proof of Theorem
\ref{thm:connect4-1}.
\hfill $\Box$\\

\begin{corollary}
\label{cor:connect4-1} With reference to Notation \ref{setup}, the
following (i), (ii) hold for $0 \leq i,j\leq D$.
\begin{enumerate}
\item ${\tilde V}^{\downarrow \downarrow}_{i,j}=0$ if $i+j<D$.
\item ${\tilde V}^{\uparrow \uparrow}_{i,j}=0$ if $i+j>D$.
\end{enumerate}
\end{corollary}

\noindent {\it Proof:} (i) This is just \cite[Theorem~6.2(ii)]{ds}.
\\
(ii) Immediate from Lemma \ref{etads} and Theorem
\ref{thm:connect4-1}(ii). \hfill
$\Box$\\

\medskip
\noindent In \cite{duality} we showed that with respect to the
standard Hermitian form the ($\downarrow,\downarrow$)--split
decomposition (resp. ($\downarrow,\uparrow$)--split decomposition)
and the ($\uparrow,\uparrow$)--split decomposition (resp.
($\uparrow,\downarrow$)--split decomposition) are dual in the
following sense.

\begin{theorem}
\label{thm:duality} {\rm \cite [Theorem~4.8]{duality}} With
reference to Notation \ref{setup} and Definition \ref{def:vtilde},
the following (i), (ii) hold for $0 \leq i,j,r,s\leq D$.
\begin{enumerate}
\item ${{\tilde V}^{\downarrow \downarrow}_{i,j}}$ and ${{\tilde V}^{\uparrow \uparrow}_{r,s}}$ are orthogonal unless $i+r=D$ and $j+s=D$.
\item ${{\tilde V}^{\downarrow \uparrow}_{i,j}}$ and ${{\tilde V}^{\uparrow \downarrow}_{r,s}}$ are orthogonal unless $i+r=D$ and $j+s=D$.
\end{enumerate}
\end{theorem}

\section{The matrices $E^{\downarrow \downarrow}_{i,j}$, $E^{\uparrow \downarrow}_{i,j}$, $E^{\downarrow \uparrow}_{i,j}$, $E^{\uparrow \uparrow}_{i,j}$}

\noindent In this section we use the split decompositions to define
the matrices $E^{\mu \nu}_{i,j}$ for $0 \leq i,j \leq D$ and $\mu,
\nu \in \lbrace \downarrow, \uparrow \rbrace$. We then discuss some
basic properties of these matrices. We start with a definition.

\begin{definition}
\label{def:Eij}
 \rm With reference to Notation \ref{setup}, for $0 \leq i,j \leq D$ and for $\mu, \nu \in \lbrace \downarrow, \uparrow
\rbrace$ we define $E^{\mu \nu}_{i,j} \in \hbox{Mat}_X(\C)$ so that
\begin{eqnarray*}
\label{Eij1} (E^{\mu \nu}_{i,j} - I){\tilde V}^{\mu \nu}_{i,j}&=&0,
\\
\label{Eij2} E^{\mu \nu}_{i,j}{\tilde V}^{\mu \nu}_{r,s}&=&0 ~
\hbox{if}~(i,j) \neq (r,s) \qquad \qquad (0 \leq r,s \leq D).
\end{eqnarray*}
In other words $E^{\mu \nu}_{i,j}$ is the projection from $V$ onto
${\tilde V}^{\mu \nu}_{i,j}$. We note that
\begin{eqnarray}
\label{EijV} E^{\mu \nu}_{i,j} V &=& {\tilde V}^{\mu \nu}_{i,j}.
\end{eqnarray}
\end{definition}

\begin{lemma}
\label{lem:prosumEij} With reference to Notation \ref{setup}, the
following (i), (ii) hold for $\mu, \nu \in \lbrace \downarrow,
\uparrow \rbrace$ and $0 \leq i,j,r,s\leq D$.
\begin{enumerate}
\item $\sum_{i=0}^D \sum_{j=0}^D E^{\mu \nu}_{i,j}=I$.
\item $E^{\mu \nu}_{i,j}E^{\mu \nu}_{r,s}=\delta_{ir} \delta_{js} E^{\mu \nu}_{i,j}$.
\end{enumerate}
\end{lemma}

\noindent {\it Proof:} Immediate from Corollary \ref{cor:drsum} and
Definition \ref{def:Eij}.
\hfill $\Box $\\

\begin{corollary}
\label{cor:sumEij1} With reference to Notation \ref{setup} and
Definition \ref{def:Eij}, the following (i), (ii) hold for $0 \leq
\eta \leq D$.
\begin{enumerate}
\item $\phi_{\eta} = \sum {E}^{\downarrow \downarrow}_{i,j}$, where the sum is
over all ordered pairs $i,j$ such that $0 \leq i,j\leq D$ and
$i+j=D+\eta$.
\item $\phi_{\eta} = \sum {E}^{\uparrow \uparrow}_{i,j}$, where
the sum is over all ordered pairs $i,j$ such that $0 \leq i,j\leq D$
and $i+j=D-\eta$.
\end{enumerate}
\end{corollary}

\noindent {\it Proof:} Combine Theorem \ref{thm:connect4-1} and
Definition \ref{def:Eij}.
\hfill $\Box $\\

\begin{corollary}
\label{cor:sumEij2} With reference to Notation \ref{setup} and
Definition \ref{def:Eij}, the following (i), (ii) hold for $-D \leq
\zeta \leq D$.
\begin{enumerate}
\item $\psi_{\zeta}= \sum {E}^{\downarrow \uparrow}_{i,j}$, where
the sum is over all ordered pairs $i,j$ such that $0 \leq i,j\leq D$
and $i+j=D+\zeta$.
\item $\psi_{\zeta}= \sum {E}^{\uparrow \downarrow}_{i,j}$, where the sum is
over all ordered pairs $i,j$ such that $0 \leq i,j\leq D$ and
$i+j=D-\zeta$.
\end{enumerate}
\end{corollary}

\noindent {\it Proof:} Combine Theorem \ref{thm:connect4-2} and
Definition \ref{def:Eij}.
\hfill $\Box $\\

\begin{corollary}
\label{cor:downupEij} With reference to Notation \ref{setup} and
Definition \ref{def:Eij}, the following (i), (ii) hold for $0 \leq
i,j\leq D$.
\begin{enumerate}
\item ${E}^{\downarrow \downarrow}_{i,j}=0$ if $i+j<D$.
\item
${E}^{\uparrow \uparrow}_{i,j}=0$ if $i+j>D$.
\end{enumerate}
\end{corollary}

\noindent {\it Proof:} Immediate from Corollary \ref{cor:connect4-1}
and Definition \ref{def:Eij}.
\hfill $\Box $\\

\begin{lemma}
\label{lem:vconj} With reference to Notation \ref{setup} and
Definition \ref{def:vtilde}, pick any $\mu, \nu \in \lbrace
\downarrow, \uparrow \rbrace$ and $0 \leq i,j\leq D$. Then the
following (i), (ii) hold for any $v \in V$.
\begin{enumerate}
\item $v \in V^{\mu \nu}_{i,j}$ if and only if ${\overline v} \in V^{\mu \nu}_{i,j}$.
\item
$v \in {\tilde V}^{\mu \nu}_{i,j}$ if and only if ${\overline v} \in
{\tilde V}^{\mu \nu}_{i,j}$.
\end{enumerate}
\end{lemma}

\noindent {\it Proof:} (i) The result follows by Definition
\ref{def:updown} and since $E_i, E^*_i$ are real for $0 \leq i \leq
D$.
\\
(ii) Routine using Definition \ref{def:vtilde} and (i) above.
 \hfill
$\Box$\\

\begin{theorem}
\label{thm:tranconjuEij} With reference to Notation \ref{setup} and
Definition \ref{def:Eij}, the following (i)--(iii) hold for $\mu,
\nu \in \lbrace \downarrow, \uparrow \rbrace$ and $0 \leq i,j\leq
D$.
\begin{enumerate}
\item $\overline{E^{\mu \nu}_{i,j}}=E^{\mu \nu}_{i,j}$.
\item $(E^{\downarrow \downarrow}_{i,j})^t=E^{\uparrow \uparrow}_{D-i,D-j}$.
\item $(E^{\downarrow \uparrow}_{i,j})^t=E^{\uparrow \downarrow}_{D-i,D-j}$.
\end{enumerate}
\end{theorem}

\noindent {\it Proof:} (i) By Definition \ref{def:Eij} it suffices
to show that $(\overline{E^{\mu \nu}_{i,j}} - I){\tilde V}^{\mu
\nu}_{i,j}=0$ and $\overline{E^{\mu \nu}_{i,j}}{\tilde V}^{\mu
\nu}_{r,s}=0$ if $(i,j) \neq (r,s)$. To see $(\overline{E^{\mu
\nu}_{i,j}} - I){\tilde V}^{\mu \nu}_{i,j}=0$, pick $v$ in ${\tilde
V}^{\mu \nu}_{i,j}$ and note that $\overline{v} \in {\tilde V}^{\mu
\nu}_{i,j}$ by Lemma \ref{lem:vconj}(ii). Now

\begin{eqnarray*}
(\overline{E^{\mu \nu}_{i,j}} - I)v &=& \overline{({E^{\mu
\nu}_{i,j}} - I)\overline{v}}\\
&=& 0
\end{eqnarray*}
so $(\overline{E^{\mu \nu}_{i,j}} - I){\tilde V}^{\mu \nu}_{i,j}=0$.
Next we fix $r,s$ ($0 \leq r,s\leq D$) such that $(r,s) \neq (i,j)$
and show that $\overline{E^{\mu \nu}_{i,j}}{\tilde V}^{\mu
\nu}_{r,s}=0$. Pick $v$ in ${\tilde V}^{\mu \nu}_{r,s}$ and note
that $\overline{v} \in {\tilde V}^{\mu \nu}_{r,s}$ by Lemma
\ref{lem:vconj}(ii). Observe that
\begin{eqnarray*}
\overline{E^{\mu \nu}_{i,j}}v &=& \overline{{E^{\mu
\nu}_{i,j}}\overline{v}}\\
&=& 0
\end{eqnarray*}
so $\overline{E^{\mu \nu}_{i,j}}{\tilde V}^{\mu \nu}_{r,s}=0$ and the result follows.\\
(ii) We show that $(E^{\uparrow \uparrow}_{D-i,D-j})^t=E^{\downarrow
\downarrow}_{i,j}$. By Definition \ref{def:Eij} it suffices to show
that $((E^{\uparrow \uparrow}_{D-i,D-j})^t - I){\tilde
V}^{\downarrow \downarrow}_{i,j}=0$ and $(E^{\uparrow
\uparrow}_{D-i,D-j})^t{\tilde V}^{\downarrow \downarrow}_{r,s}=0$ if
$(i,j) \neq (r,s)$. First we fix $r,s$ ($0 \leq r,s\leq D$) such
that $(r,s) \neq (i,j) $ and show that $(E^{\uparrow
\uparrow}_{D-i,D-j})^t{\tilde V}^{\downarrow \downarrow}_{r,s}=0$.
To do this it suffices to show that $\langle (E^{\uparrow
\uparrow}_{D-i,D-j})^t {\tilde V}^{\downarrow \downarrow}_{r,s},\, V
\rangle = 0$. Observe that
\begin{eqnarray*}
\langle (E^{\uparrow \uparrow}_{D-i,D-j})^t {\tilde V}^{\downarrow
\downarrow}_{r,s},\, V \rangle &=& \langle {\tilde V}^{\downarrow
\downarrow}_{r,s},\, E^{\uparrow \uparrow}_{D-i,D-j} V \rangle \,
\qquad \qquad \,(\hbox{by (\ref{SHF}) and (i)})\\
&=& \langle {\tilde V}^{\downarrow \downarrow}_{r,s},\, {\tilde
V}^{\uparrow \uparrow}_{D-i,D-j} \rangle \qquad \qquad \quad \,\,
(\hbox{by
(\ref{EijV})})\\
&=& 0 \, \qquad \qquad \qquad \qquad \qquad \quad \, ~(\hbox{by
Theorem \ref{thm:duality}(i)})
\end{eqnarray*}
so $(E^{\uparrow \uparrow}_{D-i,D-j})^t{\tilde V}^{\downarrow
\downarrow}_{r,s}=0$. To obtain $((E^{\uparrow
\uparrow}_{D-i,D-j})^t - I){\tilde V}^{\downarrow
\downarrow}_{i,j}=0$ combine Lemma \ref{lem:prosumEij}(i) and our
previous comments.
\\
(iii) Similar to the proof of (ii).
\hfill $\Box $\\

\begin{corollary}
\label{cor:tranconj_phi}
 With reference to Notation \ref{setup} and Definition
\ref{def:Eij}, the following (i), (ii) hold for $0 \leq \eta \leq
D$.
\begin{enumerate}
\item $\overline{\phi_{\eta}} = \phi_{\eta}$.
\item $\phi^{t}_{\eta} = \phi_{\eta}$.
\end{enumerate}
\end{corollary}

\noindent {\it Proof:} (i) Immediate from Corollary
\ref{cor:sumEij1} and Theorem \ref{thm:tranconjuEij}(i).
\\
(ii) In the equation of Corollary \ref{cor:sumEij1}(i) take the
transpose and evaluate the result using Corollary
\ref{cor:sumEij1}(ii) and Theorem \ref{thm:tranconjuEij}(ii).
\hfill $\Box $\\

\begin{corollary}
\label{cor:tranconj_psi} With reference to Notation \ref{setup} and
Definition \ref{def:Eij}, the following (i), (ii) hold for $-D \leq
\zeta \leq D$.
\begin{enumerate}
\item $\overline{\psi_{\zeta}} = \psi_{\zeta}$.
\item $\psi^{t}_{\zeta} = \psi_{\zeta}$.
\end{enumerate}
\end{corollary}

\noindent {\it Proof:} Similar to the proof of Corollary
\ref{cor:tranconj_phi}.
 \hfill $\Box $\\

\section{A restriction on the intersection numbers}
\noindent For the rest of this paper we impose the following
restriction on the intersection numbers of $\Gamma$.

\begin{assumption}
\label{def:sdcp} \rm We fix $b,\beta \in \C$ such that $b\not=1$,
and assume $\Gamma$ has classical parameters $(D,b,\alpha,\beta)$
with $\alpha=b-1$.
 This means that the intersection numbers
of $\Gamma$ satisfy
\begin{eqnarray*}
c_i &=& b^{i-1}\frac{b^i-1}{b-1}, \\
b_i &=& (\beta+1-b^i)\frac{b^D-b^i}{b-1}
\end{eqnarray*}
for $0 \leq i \leq D$ \cite[p.~193]{bcn}. We remark that $b$ is an
integer and $b\not=0$, $b\not=-1$ \cite[Proposition~6.2.1]{bcn}. For
notational convenience we fix $q \in \C$ such that $b = q^2$. We
note that $q$ is nonzero and not a root of unity.
\end{assumption}

\begin{remark}
\label{rem:qp} \rm Referring to Assumption \ref{def:sdcp}, the
restriction $\alpha=b-1$
 implies that
$\Gamma$ is formally self-dual \cite[Corollary~8.4.4]{bcn}.
Consequently there exists an ordering $E_0, E_1, \ldots, E_D$ of the
primitive idempotents of $\Gamma$, with respect to which the Krein
parameter $q^h_{ij}$ is equal to the intersection number $p^h_{ij}$
for $0 \leq h,i,j\leq D$. In particular $\Gamma$ is $Q$-polynomial
with respect to $E_0, E_1, \ldots, E_D$. We fix this ordering of the
primitive idempotents for the rest of the paper.
\end{remark}

\medskip
\noindent With reference to Assumption \ref{def:sdcp}, Ito and
Terwilliger displayed an action of the $q$-tetrahedron algebra
$\boxtimes_q$ on the standard module of $\Gamma$ \cite{drgqtet}. To
describe this action they defined eight matrices in
 $\hbox{Mat}_X(\C)$, called
\begin{eqnarray}
\label{eq:list} A,\quad A^*,\quad B,\quad B^*, \quad K,\quad
K^*,\quad \Phi,\quad \Psi.
\end{eqnarray}
For each matrix in (\ref{eq:list}) we compute the transpose and
complex conjugate. Using this infor-\\mation we compute the
transpose and complex conjugate for each generator of $\boxtimes_q$.

\section{The matrices $A$, $A^*$, $B$, $B^*$, $K$, $K^*$, $\Phi$, $\Psi$}

\noindent In this section we define the matrices (\ref{eq:list}) and
discuss their properties. We begin with a lemma.

\begin{lemma}
{\rm \cite[Corollary~8.4.4]{bcn}} \label{lem:alpha} With reference
to Assumption \ref{def:sdcp}, there exist $\alpha_0,\alpha_1 \in \C$
such that each of $\theta_i$, $\theta^*_i$ is
 $\alpha_0+\alpha_1 q^{D-2i}$ for $0 \leq i \leq D$.
Moreover $\alpha_1\not=0$.
\end{lemma}

\begin{definition}
\label{def:aas} {\rm \cite [Definition 9.2]{drgqtet}}
 \rm With reference to Assumption \ref{def:sdcp}, we define $A,A^* \in
\hbox{Mat}_X(\C)$ so that
\begin{eqnarray*}
A_1&=&\alpha_0I+\alpha_1A,
\\
A^*_1&=&\alpha_0I+\alpha_1A^*,
\end{eqnarray*}
where $\alpha_0,\alpha_1$ are from Lemma \ref{lem:alpha}. Thus for
$0 \leq i \leq D$ the space $E_iV$ (resp. $E^*_iV$) is an eigenspace
of $A$ (resp. $A^*$) with eigenvalue
 $q^{D-2i}$.
\end{definition}

\medskip
\noindent The following result is immediate from Lemma
\ref{lem:alpha} and Definition \ref{def:aas}.

\begin{lemma}
\label{lem:sumAA^*} With reference to Assumption \ref{def:sdcp} and
Definition \ref{def:aas}, the following (i), (ii) hold.
\begin{enumerate}
\item $A=\sum_{i=0}^D q^{D-2i} E_i$.
\item $A^*=\sum_{i=0}^D q^{D-2i} E^*_i$.
\end{enumerate}
\end{lemma}

\begin{definition}
{\rm \cite[Definition~10.4]{drgqtet}} \label{def:listdef} \rm With
reference to Definition \ref{def:vtilde} and Assumption
\ref{def:sdcp}, we define $B$, $B^*$, $K$, $K^*$, $\Phi$, $\Psi$ to
be the unique matrices in $\hbox{Mat}_X(\C)$ that satisfy the
requirements of the following table for $0 \leq i,j\leq D$.

\medskip
\centerline{
\begin{tabular}[t]{c|c}
        {\rm The matrix}
     &
  {\rm is $0$ on}
\\      \hline
         $B-q^{i-j}I$
         &
  ${\tilde V}_{i,j}^{\downarrow \uparrow}$
\\
 $B^*-q^{j-i}I$
         &
 ${\tilde V}_{i,j}^{\uparrow \downarrow}$
\\
 $K-q^{i-j}I$
         &
 ${\tilde V}_{i,j}^{\downarrow \downarrow}$
\\
 $K^*-q^{i-j}I$
         &
 ${\tilde V}_{i,j}^{\uparrow \uparrow}$
\\
$\Phi-q^{i+j-D}I$
         &
 ${\tilde V}_{i,j}^{\downarrow \downarrow}$
\\
 $\Psi-q^{i+j-D}I$
&
 ${\tilde V}_{i,j}^{\downarrow \uparrow}$
\end{tabular}}

\medskip
\end{definition}

\begin{proposition}
{\rm \cite[Lemma~12.1]{drgqtet}} \label{prop:central} With reference
to Assumption \ref{def:sdcp}, the following (i), (ii) hold.
\begin{enumerate}
\item Each of the matrices from the list (\ref{eq:list}) is
contained in $T$. \item Each of $\Phi, \Psi$ is central in $T$.
\end{enumerate}
\end{proposition}

\begin{proposition}
\label{prop:newlistdef} With reference to Assumption \ref{def:sdcp}
and Definition \ref{def:listdef}, the following (i)--(vi) hold.
\begin{enumerate}
\item $B = \sum_{i=0}^D \sum_{j=0}^D q^{i-j} E^{\downarrow \uparrow}_{i,j}$.
\item $B^* = \sum_{i=0}^D \sum_{j=0}^D q^{j-i} E^{\uparrow \downarrow}_{i,j}$.
\item $K = \sum_{i=0}^D \sum_{j=0}^D q^{i-j} E^{\downarrow \downarrow}_{i,j}$.
\item $K^* = \sum_{i=0}^D \sum_{j=0}^D q^{i-j} E^{\uparrow \uparrow}_{i,j}$.
\item $\Phi = \sum_{i=0}^D \sum_{j=0}^D q^{i+j-D} E^{\downarrow \downarrow}_{i,j}$.
\item $\Psi = \sum_{i=0}^D \sum_{j=0}^D q^{i+j-D} E^{\downarrow \uparrow}_{i,j}$.
\end{enumerate}
\end{proposition}

\noindent {\it Proof:} (i) Let $\check{B}=\sum_{i=0}^D \sum_{j=0}^D
q^{i-j} E^{\downarrow \uparrow}_{i,j}$. We show that $B=\check{B}$.
Using Definition \ref{def:Eij} we find that $\check{B}-q^{i-j}I$ is
0 on ${\tilde V}_{i,j}^{\downarrow \uparrow}$ for $0 \leq i,j\leq
D$, so $B=\check{B}$ in view of Definition \ref{def:listdef}.\\
(ii)--(vi) Similar to the proof of (i).
\hfill $\Box $\\

\begin{lemma}
\label{lem:symPPetazeta} With reference to Assumption \ref{def:sdcp}
and Definition \ref{def:listdef}, the following (i)--(iv) hold.
\begin{enumerate}
\item $\Phi = \sum_{\eta=0}^D  q^{\eta} \phi_\eta$.
\item $\Phi = \sum_{i=0}^D \sum_{j=0}^D q^{D-i-j} E^{\uparrow \uparrow}_{i,j}$.
\item $\Psi = \sum_{\zeta=-D}^D  q^{\zeta} \psi_\zeta$.
\item $\Psi = \sum_{i=0}^D \sum_{j=0}^D q^{D-i-j} E^{\uparrow \downarrow}_{i,j}$.
\end{enumerate}
\end{lemma}

\noindent {\it Proof:} (i) Evaluate Proposition
\ref{prop:newlistdef}(v) using Corollary \ref{cor:sumEij1}(i).
\\
(ii) Evaluate (i) above using Corollary \ref{cor:sumEij1}(ii).
\\
(iii) Evaluate Proposition \ref{prop:newlistdef}(vi) using Corollary
\ref{cor:sumEij2}(i).
\\
(iv) Evaluate (iii) above using Corollary \ref{cor:sumEij2}(ii).
\hfill $\Box $\\

\section{The transpose of $A$, $A^*$, $B$, $B^*$, $K$, $K^*$, $\Phi$, $\Psi$}
In this section we find the transpose of each matrix in the list
(\ref{eq:list}). We start with $A$ and $A^*$.

\begin{theorem}
\label{thm:symA} With reference to Assumption \ref{def:sdcp} and
Definition \ref{def:aas}, the following (i), (ii) hold.
\begin{enumerate}
\item $A$ is symmetric.
\item $A^*$ is symmetric.
\end{enumerate}
\end{theorem}

\noindent {\it Proof:} Immediate from Definition \ref{def:aas} and
since each of $A_1$, $A^*_1$ is symmetric.
\hfill $\Box $\\

\begin{theorem}
\label{thm:tranBK} With reference to Assumption \ref{def:sdcp} and
Definition \ref{def:listdef}, the following (i), (ii) hold.
\begin{enumerate}
\item $B^t=B^*$.
\item $K^t={K^*}^{-1}$.
\end{enumerate}
\end{theorem}

\noindent {\it Proof:} (i) Combining Theorem
\ref{thm:tranconjuEij}(iii) and Proposition
\ref{prop:newlistdef}(i),(ii) we have
\begin{eqnarray*}
B^t &=& \sum_{i=0}^D \sum_{j=0}^D q^{i-j} E^{\uparrow \downarrow}_{D-i,D-j}\\
        &=& \sum_{i=0}^D \sum_{j=0}^D q^{(D-j)-(D-i)} E^{\uparrow \downarrow}_{D-i,D-j}\\
        &=& \sum_{i=0}^D \sum_{j=0}^D q^{j-i} E^{\uparrow \downarrow}_{i,j}\\
        &=& B^*.
\end{eqnarray*}
\\
(ii) Using Lemma \ref{lem:prosumEij} and Proposition
\ref{prop:newlistdef}(iv) we find ${K^*}^{-1}=\sum_{i=0}^D
\sum_{j=0}^D q^{j-i} E^{\uparrow \uparrow}_{i,j}$. Combining this
with Theorem \ref{thm:tranconjuEij}(ii) and Proposition
\ref{prop:newlistdef}(iii) we obtain
\begin{eqnarray*}
K^t &=& \sum_{i=0}^D \sum_{j=0}^D q^{i-j} E^{\uparrow \uparrow}_{D-i,D-j}\\
        &=& \sum_{i=0}^D \sum_{j=0}^D q^{(D-j)-(D-i)} E^{\uparrow \uparrow}_{D-i,D-j}\\
        &=& \sum_{i=0}^D \sum_{j=0}^D q^{j-i} E^{\uparrow \uparrow}_{i,j}\\
        &=& {K^*}^{-1}.
\end{eqnarray*}
\hfill
$\Box $\\

\begin{theorem}
\label{thm:symPP} With reference to Assumption \ref{def:sdcp} and
Definition \ref{def:listdef}, the following (i), (ii) hold.
\begin{enumerate}
\item $\Phi$ is symmetric.
\item $\Psi$ is symmetric.
\end{enumerate}
\end{theorem}

\noindent {\it Proof:} (i) Immediate from Corollary
\ref{cor:tranconj_phi}(ii) and Lemma \ref{lem:symPPetazeta}(i).
\\
(ii) Immediate from Corollary \ref{cor:tranconj_psi}(ii) and Lemma
\ref{lem:symPPetazeta}(iii).
\hfill $\Box $\\

\medskip
\noindent We finish this section with a comment.

\begin{proposition}
\label{prop:invPP} With reference to Definition \ref{def:Eij} and
Assumption \ref{def:sdcp}, the following (i), (ii) hold.
\begin{enumerate}
\item $\Phi^{-1}=\sum_{i=0}^D \sum_{j=0}^D q^{i+j-D} E^{\uparrow \uparrow}_{i,j}$.
\item $\Psi^{-1}=\sum_{i=0}^D \sum_{j=0}^D q^{i+j-D} E^{\uparrow \downarrow}_{i,j}$.
\end{enumerate}
\end{proposition}

\noindent {\it Proof:} (i) Define $\Lambda=\sum_{i=0}^D \sum_{j=0}^D
q^{i+j-D} E^{\uparrow \uparrow}_{i,j}$. Using Lemma
\ref{lem:prosumEij} and Lemma \ref{lem:symPPetazeta}(ii) we
routinely find $\Phi \Lambda = I$, so $\Phi^{-1}=\Lambda$.
\\
(ii) Similar to the proof of (i).
\hfill $\Box $\\

\section{The complex conjugate of $A$, $A^*$, $B$, $B^*$, $K$, $K^*$, $\Phi$, $\Psi$}

\noindent In this section we find the complex conjugate of each
matrix in the list (\ref{eq:list}). We begin with a remark.

\begin{remark}
\label{rem:qp} \rm Recall from Assumption \ref{def:sdcp} that $q \in
\C$ satisfies $b = q^2$. Define $q{'}=-q$ and note that $(q')^2=b$.
There are two cases to consider. For $b>1$ we have $\overline{q}=q$
and $\overline{q'}=q'$. For $b<-1$ we have $q \in i\mathbb R$ so
$\overline{q}=q'$. For each matrix $S$ obtained using $q$ let $S'$
denote the corresponding matrix associated with $q'$.
\end{remark}

\begin{lemma}
\label{lem:2roots} With reference to Assumption \ref{def:sdcp} and
Definition \ref{def:aas}, the following (i)--(iv) hold.
\begin{enumerate}
\item If $D$ is even then $A=A'$.
\item If $D$ is odd then $A=-A'$.
\item If $D$ is even then $A^*={A^*}'$.
\item If $D$ is odd then $A^*=-{A^*}'$.
\end{enumerate}
\end{lemma}

\noindent {\it Proof:} (i), (ii) By Lemma \ref{lem:sumAA^*}(i) and
Remark \ref{rem:qp} we find
\begin{eqnarray*}
A' &=& \sum_{i=0}^D (q')^{D-2i} E_i
\\
&=& \sum_{i=0}^D (-q)^{D-2i} E_i.
\end{eqnarray*}
Comparing this with Lemma \ref{lem:sumAA^*}(i) we get the result.\\
(iii), (iv) Similar to the proof of (i), (ii) above.
\hfill $\Box $\\

\begin{theorem}
\label{thm:realimagconj} With reference to Assumption \ref{def:sdcp}
and Definitions \ref{def:aas}, \ref{def:listdef}, the following (i),
(ii) hold.
\begin{enumerate}
\item Assume $b>1$. Then each of $A$, $A^*$, $B$, $B^*$, $K$, $K^*$, $\Phi$, $\Psi$
is real.
\item Assume $b<-1$. Then
\begin{align}
\label{imagconj1}
 \overline{A}&=A',& &\overline{A^*}={A^*}', \\
\label{imagconj2}
 \overline{B}&=B',& &\overline{B^*}={B^*}',
\\
\label{imagconj3}
 \overline{K}&=K',& &\overline{K^*}={K^*}',
\\
\label{imagconj4}
 \overline{\Phi}&=\Phi',& &\,\,\,\overline{\Psi}=\Psi'.
\end{align}
\end{enumerate}
\end{theorem}

\noindent {\it Proof:} (i) Each of $A$, $A^*$ is real by Lemma
\ref{lem:sumAA^*}, since $q$ is real and each of $E_i$, $E^*_i$ is
real for $0 \leq i \leq D$. For the remaining matrices the results
are immediate from Theorem \ref{thm:tranconjuEij}(i), Proposition
\ref{prop:newlistdef}, and Remark \ref{rem:qp}.
\\
(ii) We first obtain the equation on the left in (\ref{imagconj1}).
By Lemma \ref{lem:sumAA^*}(i) and Remark \ref{rem:qp} we find
\begin{eqnarray*}
A' &=& \sum_{i=0}^D (q')^{D-2i} E_i
\\
&=& \sum_{i=0}^D (\overline{q})^{D-2i} E_i
\end{eqnarray*}
and this equals $\overline{A}$ since $E_i$ is real for $0 \leq i
\leq D$. We have now obtained the equation on the left in
(\ref{imagconj1}). The equation on the right in (\ref{imagconj1}) is
similarly obtained. Next we obtain the equation on the left in
(\ref{imagconj2}). By Theorem \ref{thm:tranconjuEij}(i), Proposition
\ref{prop:newlistdef}(i), and Remark \ref{rem:qp} we have
\begin{eqnarray*}
B'&=& \sum_{i=0}^D \sum_{j=0}^D (q')^{i-j} E^{\downarrow \uparrow}_{i,j}\\
                  &=&\sum_{i=0}^D \sum_{j=0}^D (\overline{q})^{i-j} E^{\downarrow \uparrow}_{i,j}\\
                  &=&\overline{B}.
\end{eqnarray*}
The remaining equations are similarly obtained.
\hfill $\Box $\\

\section{The $q$-tetrahedron algebra $\boxtimes_q$}
\noindent Referring to Assumption \ref{def:sdcp}, Ito and
Terwilliger displayed an action of the $q$-tetrahedron algebra
$\boxtimes_q$ on the standard module $V$ of $\Gamma$ \cite{drgqtet}.
In this section we compute the tranpose and complex conjugate for
the action of each generator of $\boxtimes_q$ on $V$. We start with
the definition of $\boxtimes_q$.

\medskip
\noindent For a nonzero scalar $q \in \C$ such that $q^2\not=1$ we
define
\begin{eqnarray*}
\lbrack n \rbrack_q = \frac{q^n-q^{-n}}{q-q^{-1}}, \qquad \qquad n =
0,1,2,\ldots \label{eq:nbrack}
\end{eqnarray*}
We let $\Z_4 = \Z/4\Z$ denote the cyclic group of order 4.

\begin{definition}
\rm \cite[Definition 10.1]{qtet} \label{def:qtet} Let $\boxtimes_q$
denote the unital associative $\C$-algebra that has generators
\begin{eqnarray*}
\lbrace x_{ij}\;|\; i,j \in \Z_4,\;j-i=1 \;\mbox{or} \;j-i=2\rbrace
\end{eqnarray*}
and the following relations:
\begin{enumerate}
\item For $i,j\in \Z_4$ such that $j-i=2$,
\begin{eqnarray*}
x_{ij}x_{ji} = 1. \label{eq:qrel0}
\end{eqnarray*}
\item For $h,i,j\in \Z_4$ such that the pair $(i-h,j-i)$ is one of
$(1,1), (1,2), (2,1)$,
\begin{eqnarray*}
\frac{qx_{hi}x_{ij}-q^{-1}x_{ij}x_{hi}}{q-q^{-1}}=1.
\label{eq:qrel1}
\end{eqnarray*}
\item For $h,i,j,k\in \Z_4$ such that $i-h=j-i=k-j=1$,
\begin{eqnarray*}
\label{eq:qserre} x_{hi}^3x_{jk} - \lbrack 3 \rbrack_q
x_{hi}^2x_{jk}x_{hi} + \lbrack 3 \rbrack_q x_{hi}x_{jk}x_{hi}^2-
x_{jk}x_{hi}^3=0.
\end{eqnarray*}
\end{enumerate}
We call $\boxtimes_q$ the {\it $q$-tetrahedron algebra} or
``$q$-tet'' for short.
\end{definition}

\begin{proposition}
{\rm \cite[Theorem 11.1]{drgqtet}} \label{prop:mres} With reference
to Assumption \ref{def:sdcp}, there exists a $\boxtimes_q$-module
structure on $V$ such that the generators $x_{ij}$ act as follows:

\medskip
\centerline{
\begin{tabular}[t]{c|cccccccc}
        {\rm generator}
        & $x_{01}$
         & $x_{12}$
         & $x_{23}$
         & $x_{30}$
         & $x_{02}$
         & $x_{20}$
         & $x_{13}$
         & $x_{31}$
    \\
    \hline
{\rm action on $V$} &  $A\Phi \Psi^{-1}$ & $B \Phi^{-1}$ & $A^*
\Phi\Psi$ & $B^*\Phi^{-1}$ & $K\Psi^{-1}$ &  $\Psi K^{-1}$
&$K^*\Psi$ &$\Psi^{-1}{K^*}^{-1}$
\end{tabular}}
\end{proposition}

\begin{theorem}
\label{thm:tranx_ij} With reference to Assumption \ref{def:sdcp} and
Proposition \ref{prop:mres}, the following hold on $V$.
\begin{eqnarray*}
x^{t}_{01}=x_{01}, \qquad x^{t}_{12}=x_{30}, \qquad
x^{t}_{23}=x_{23}, \qquad x^{t}_{30}=x_{12},\\
x^{t}_{02}=x_{31}, \qquad x^{t}_{20}=x_{13}, \qquad
x^{t}_{13}=x_{20},\qquad x^{t}_{31}=x_{02}.
\end{eqnarray*}
\end{theorem}

\noindent {\it Proof:} Combine Proposition \ref{prop:central},
Theorems \ref{thm:symA}--\ref{thm:symPP}, and Proposition
\ref{prop:mres}.
\hfill $\Box $\\

\begin{theorem}
With reference to Assumption \ref{def:sdcp} and Proposition
\ref{prop:mres}, the following (i), (ii) hold on $V$.
\begin{enumerate}
\item Assume $b>1$. Then $\overline{x_{ij}}=x_{ij}$ for each generator $x_{ij}$ of $\boxtimes_q$.
\item Assume $b<-1$. Then $\overline{x_{ij}}=x'_{ij}$ for each
generator $x_{ij}$ of $\boxtimes_q$.
\end{enumerate}
\end{theorem}

\noindent {\it Proof:} (i) Immediate from Theorem
\ref{thm:realimagconj}(i) and Proposition \ref{prop:mres}.
\\
(ii) Immediate from Theorem \ref{thm:realimagconj}(ii) and
Proposition \ref{prop:mres}.
\hfill $\Box $\\

\bigskip
\noindent{\Large\bf Acknowledgements}

\medskip
\noindent The author would like to thank Professor Paul M.
Terwilliger for his valuable comments.

\noindent Joohyung Kim \\
National Institute for Mathematical Sciences\\
385-16, Doryong-dong, Yuseong-gu,\\
Daejeon 305-340, Republic of Korea\\
email: jayhkim@nims.re.kr

\end{document}